
\input amstex
\documentstyle{amsppt}
\magnification=\magstep1
\NoBlackBoxes

\def\sgn{\text{sgn}}
\topmatter

\title Positive definite distributions  
and subspaces of $L_{-p}$ with applications 
to stable processes
\endtitle

\author Alexander Koldobsky \endauthor 
\address Division of Mathematics and Statistics, 
University
of Texas at San Antonio, San Antonio, TX 78249, U.S.A. \endaddress
\email koldobsk\@ringer.cs.utsa.edu \endemail

\rightheadtext{Positive definite distributions}

\abstract  We define embedding of an $n$-dimensional 
normed space into 
$L_{-p},\ 0<p<n$ by extending analytically with respect to $p$
the corresponding property of the classical $L_p$-spaces.
The well-known connection between embeddings into
$L_p$ and positive definite functions is extended to 
the case of negative $p$ by showing that a normed 
space embeds in $L_{-p}$ if and only if $\|x\|^{-p}$ is a 
positive definite distribution.
Using this criterion, we generalize the recent solutions
to the 1938 Schoenberg's problems by proving that the spaces $\ell_q^n,
\ 2<q\le \infty$ embed in $L_{-p}$ if and only if $p\in [n-3,n).$
We show that the technique of embedding in $L_{-p}$ can 
be applied to stable processes in some situations where 
standard methods do not work. As an example, we prove inequalities
of correlation type for the expectations of norms of stable vectors.
In particular, for every $p\in [n-3,n),$
$\Bbb E(\max_{i=1,...,n} |X_i|^{-p}) \ge 
\Bbb E(\max_{i=1,...,n} |Y_i|^{-p}),$
where $X_1,...,X_n$ and $Y_1,...,Y_n$ are 
jointly $q$-stable symmetric random variables, $0<q\le 2,$ 
so that, for some $k\in \Bbb N,\ 1\le k <n,$ 
the vectors $(X_1,...,X_k)$ and $(X_{k+1},...,X_n)$ 
have the same distributions as $(Y_1,...,Y_k)$ and 
$(Y_{k+1},...,Y_n),$ respectively, but $Y_i$ and $Y_j$ are
independent for every choice of $1\le i\le k,\ k+1\le j\le n.$
\endabstract


\thanks Research supported in part 
by the NSF Grant DMS-9531594 \endthanks

\endtopmatter \document \baselineskip=14pt

\head 1. Introduction \endhead

The connections between stable measures, positive definite norm
dependent functions and embedding of normed spaces in $L_p$ 
were discovered by P.Levy \cite {12} as parts of his theory
of stable processes, and, since then, those connections
have been under intensive development (see \cite{10,15} for
the most recent surveys). 
In particular, P.Levy pointed out that an 
$n$-dimensional normed space $B=(\Bbb R^n,\|\cdot\|)$ 
embeds isometrically in $L_p,\ p>0$ if and only if there
exists a finite Borel measure $\gamma$ on the unit sphere
$\Omega$ in $\Bbb R^n$ so that
$$\|x\|^p =\int_S |(x,\xi)|^p\ d\gamma(\xi) \tag{1}$$
for every $x\in \Bbb R^n.$ On the other hand, for $0<p\le 2,$
the representation (1) exists if and only if the function
$\exp(-\|x\|^p)$ is positive definite and, hence, 
is the characteristic function of a symmetric stable measure 
in $\Bbb R^n.$ We call (1) the Blaschke-Levy 
representation of the norm with the exponent $p$ and
measure $\gamma$ (see \cite{11} for the history, generalizations 
and applications of this representation). 

Several applications of the Blaschke-Levy
representation to stable processes depend on the standard
procedure of using (1) to estimate the expectation
of the norm of a stable vector (we give an example in Section 4). 
Usually, those applications do not use the Banach 
space structure of the space $L_p,$ and they work equally
well for $p\ge 1$ and $p\in (0,1).$ Moreover, 
when $p<2$ becomes smaller one can expect more normed spaces
to admit the representation (1) with the exponent $p,$ because
for $0<p_1<p_2\le 2,$ the space $L_{p_2}$ embeds isometrically 
in $L_{p_1}$ (see \cite{1}).
However, the spaces $\ell_\infty^n,
\ n\ge 3$ do not embed in any of the spaces $L_p$ with $p>0,$ 
and the spaces $\ell_q^n,\ n\ge 3,\ q>2$ do not embed in $L_p$
with $0<p\le 2$ (see \cite{14,7}; note that the latter results 
solved the 1938 Schoenberg's problems on positive definite
functions \cite{19}.) These spaces (especially $\ell_\infty^n$) 
are particularly important in the theory of stable processes,
and it seems to be natural to try to modify the standard
technique so that it works for those spaces.

These were the reasons which led the author to an attempt 
to get more norms involved by generalizing the Blaschke-Levy 
representation (and embedding in $L_p$) to the case of
negative $p.$ In Section 2, we define the
Blaschke-Levy representation in $\Bbb R^n$ with negative 
exponents $-p,\ 0<p<n,$ and we say that the existence of 
such a representation for a normed space means that the space 
embeds in $L_{-p}.$ The definition is "analytic" with respect 
to $p,$ which might allow us to transfer properties of the 
spaces $L_p$ in both directions between the
positive and negative values of $p.$ We show that the connection
between embeddings in $L_p$ and positive definiteness remains
in force, namely, a space $B=(\Bbb R^n,\|\cdot\|)$ embeds in 
$L_{-p}$ if and only if $\|x\|^{-p}$ is a positive definite 
distribution on $\Bbb R^n.$ Recall that in the positive 
case the condition is that the distribution $\Gamma(-p/2)\|x\|^p$ 
must be positive definite outside of the origin (see \cite{8}; 
$p$ is not an even integer).

In Section 3, we show that, for $2<q\le \infty,\ 0<p<n,\ n\ge 3$ 
the function $\|x\|_q^{-p}$ is a positive definite distribution 
if and only if $p\in [n-3,n),$ where $\|x\|_q$ stands for the norm
of the space $\ell_q^n.$ This means that the spaces $\ell_\infty^n$
and $\ell_q^n,\ q>2$ embed in $L_{-p},\ 0<p<n$ if and only if 
$p\in [n-3,n),$
and this also includes the case $n=2.$ (Use the well-known fact
\cite{3,6,13}
that every two-dimensional Banach space embeds in $L_p$ for every 
$p\in (0,1]$, and then apply Theorem 2 from this
paper.)

In Section 4, we give an example of how the standard technique 
of the theory of stable processes can be modified by using 
embeddings in $L_{-p}.$ For $B=(\Bbb R^n,\|\cdot\|),\ p\in \Bbb R,$
we consider the problem of
optimization of the expectation $\Bbb E(\|X\|^p)$ of
the norm of a symmetric $q$-stable random vector $X$
in $\Bbb R^n$ in the following sense. 
Let $1\le k < n,\ 0<q\le 2$ and $X_1,...,X_n$ and $Y_1,...,Y_n$ be 
jointly $q$-stable symmetric random variables, so that 
the vectors $(X_1,...,X_k)$ and $(X_{k+1},...,X_n)$ 
have the same distributions as $(Y_1,...,Y_k)$ and 
$(Y_{k+1},...,Y_n),$ respectively, but $Y_i$ and $Y_j$ are
independent for every choice of $1\le i\le k,\ k+1\le j\le n.$
We compare the expectations $\Bbb E(\|X\|^p)$ and
$\Bbb E(\|Y\|^p)$. First, we apply the standard methods
to the case where $p>0$ and $B$ is a subspace of $L_p,$
and we prove that $\Bbb E(\|X\|^p)\le \Bbb E(\|Y\|^p)$ 
for each $p<q.$ Then, we show that the technique of embedding in
$L_{-p}$ leads to similar results for a larger class of spaces $B.$
In particular, for every $p\in [n-3,n),$
$$\Bbb E(\max_{i=1,...,k} |X_i|^{-p}) \ge 
\Bbb E(\max_{i=1,...,k} |Y_i|^{-p}).$$
The question of what happens to the latter inequality when 
the exponent $-p$ is replaced by 1 is open, and, in the Gaussian case,
this question is the matter of the weak version of the 
well-known Gaussian correlation problem (see \cite{18} for the
most recent developments).

\head 2. Positive definite distributions and embeddings in $L_{-p}$
\endhead

The main tool of this paper is the Fourier transform
of distributions. As usual, we denote by $\Cal S(\Bbb R^n)$ 
the space of rapidly decreasing infinitely differentiable 
functions (test functions) in $\Bbb R^n,$ and 
$\Cal S^{'}(\Bbb R^n)$ is
the space of distributions over $\Cal S(\Bbb R^n).$  
The Fourier transform of a distribution $f\in \Cal S^{'}(\Bbb R^n)$ 
is defined by $\langle \hat{f},\hat{\phi} \rangle = 
(2\pi)^n \langle f,\phi \rangle$ for every
test function $\phi.$  A distribution is called even homogeneous
of degree $p\in \Bbb R$  if 
$\langle f(x), \phi(x/\alpha) \rangle = |\alpha|^{n+p} 
\langle f,\phi \rangle$
for every test function $\phi$ and every 
$\alpha\in \Bbb R,\ \alpha\neq 0.$ The Fourier transform of 
an even homogeneous distribution of degree $p$ is an even 
homogeneous distribution of degree $-n-p.$
If $p>-1$ and $p$ is not an even integer, then 
the Fourier transform of the function
$h(z)=|z|^p,\ z\in \Bbb R$ is equal to 
$(|z|^p)^{\wedge}(t) = c_p |t|^{-1-p}$ 
(see \cite{4, p.173}), where
$c_p = {{2^{p+1}\sqrt{\pi}\ \Gamma((p+1)/2)}\over{\Gamma(-p/2)}}.$
The well-known connection between the Radon transform and
the Fourier transform is that, for every $\xi\in \Omega,$
the function $t\to \hat\phi(t\xi)$ 
is the Fourier transform of the function 
$z\to R\phi(\xi;z)=\int_{(x, \xi)=z} \phi (x)\, dx$
($R$ stands for the Radon transform). 
A distribution $f$ is called 
positive definite if, for every test function $\phi,$ 
$\langle f, \phi * \overline{\phi(-x)} \rangle \ge 0.$
A distribution is positive definite
if and only if it is the Fourier transform of 
a tempered measure in $\Bbb R^n$ (\cite{5, p.152}). Recall that
a (non-negative, not necessarily finite) measure $\mu$ is  
called tempered if 
$$\int_{\Bbb R^n} (1+\|x\|_2)^{-\beta}\ d\mu(x)< \infty$$
for some $\beta >0.$
Every positive distribution (in the sense that
$\langle f, \phi \rangle \ge 0$ for every non-negative
test function $\phi$) is a tempered measure \cite{5, p.147}.

\medbreak 

Throughout the paper $\|x\|$ stands for a  
homogeneous of degree 1, continuous, 
positive outside of the origin function on $\Bbb R^n.$ 
We say that $B=(\Bbb R^n, \|\cdot\|)$ is a homogeneous
$n$-dimensional space. Clearly, the class of homogeneous
spaces contains all finite dimensional normed and quasi-normed
spaces. It is easily seen that every functional $\|x\|$ is equivalent
to the Euclidean norm in the sense that, for every $x\in \Bbb R^n,$ 
$K_1\|x\|_2 \le \|x\| \le K_2 \|x\|_2$ for some positive 
constants $K_1, K_2.$ Hence, $\|x\|^{-p}$ is a locally 
integrable function on $\Bbb R^n$ for every $p\in (0,n).$

\bigbreak

Now we are ready to define the Blaschke-Levy representation
with negative exponents $p.$ 
Indeed, the formula (1) does not make sense if $p<-1.$
However, let us start with positive $p$ and apply functions 
in both sides (1) to a test function $\phi:$
$$\int_{\Bbb R^n} \|x\|^p \phi(x)\ dx =
\int_\Omega d\gamma (\xi) \int_{\Bbb R^n} |(x,\xi)|^p \phi(x)\ dx =$$
$$\int_\Omega d\gamma (\xi) \int_{\Bbb R} 
|z|^p \big(\int_{(x, \xi)=z} \phi (x)\, dx\big)\ dz =
\int_\Omega  \langle |z|^p, R\phi(\xi;z) \rangle\  d\gamma (\xi) =$$
$$c_p\int_\Omega  \langle |t|^{-1-p}, \hat\phi(t\xi) 
\rangle\ d\gamma (\xi).$$
If $p$ is negative the function $|t|^{-1-p}$ is locally integrable,
which allows to write $\langle |t|^{-1-p}, \hat\phi(t\xi) \rangle$  
as an integral, 
and this is how we extend the Blaschke-Levy representation:

\proclaim{Definition} Let $B=(\Bbb R^n,\|\cdot\|)$
be an $n$-dimensional homogeneous space, $p\in (0,n).$ 
We say that the norm of $B$ admits the Blaschhke-Levy
representation with the exponent $-p,$ if there exists a 
finite symmetric measure $\gamma$ on the sphere $\Omega$ so that, 
for every test function $\phi,$
$$\int_{\Bbb R^n} \|x\|^{-p}\phi(x)\ dx =
\int_\Omega d\gamma (\xi) \int_{\Bbb R} |t|^{p-1} \hat\phi(t\xi)\ dt.
\tag{2}$$ 
If the norm of $B$ satisfies (2) with a measure $\gamma$,
we also say that the space $B$ embeds in $L_{-p}.$
\endproclaim

\bigbreak

It is easy to show the uniqueness of the representation (2). 
In fact, consider the test functions $\phi$ of the form
$$\phi(x) = h(t) u(\xi),\ x=t\xi,\ t\in \Bbb R,\ t>0,
\ \xi\in \Omega,\tag{3}$$
where $h$ is a non-negative test function on $\Bbb R,$ and
$u$ is an infinitely differentiable even function on the sphere $\Omega.$
If a norm admits the representation (2) with two measures
$\gamma_1$ and $\gamma_2,$ then applying (2) to the test functions
whose Fourier transforms have the form (3), we get that, for every $u,$  
$$\int_\Omega u(\xi)\ d\gamma_1 (\xi) =
\int_\Omega u(\xi)\ d\gamma_2 (\xi),$$
which implies $\gamma_1 = \gamma_2.$
\bigbreak
Similar to the positive case, embedding into $L_{-p}$
is closely related to positive definiteness.
The following fact will serve as a tool for checking
whether certain spaces embed in $L_{-p}.$

\proclaim{Theorem 1} An $n$-dimensional homogeneous
space $B=(\Bbb R^n,\|\cdot\|)$ embeds in $L_{-p},\ p\in (0,n)$
if and only if $\|x\|^{-p}$ is a positive definite distribution.
\endproclaim
                                 
\demo{Proof} Suppose that $B$ embeds in $L_{-p}.$ 
For every non-negative test function $\phi,$ using (2) and the fact
that $(\hat\phi)^{\wedge}(x) = (2\pi)^n \phi(-x),$ we get
$$ \langle (\|x\|^{-p})^{\wedge}, \phi \rangle =
\int_{\Bbb R^n} \|x\|^{-p}\hat\phi(x)\ dx =  (2\pi)^n
\int_\Omega d\gamma (\xi) \int_{\Bbb R} |t|^{p-1} \phi(t\xi)\ dt \ge 0,
$$ 
which shows that $(\|x\|^{-p})^{\wedge}$ is a positive distribution 
over $\Cal S(\Bbb R^n,)$ and, hence, the distribution  
$\|x\|^{-p}$ is positive definite.

Conversely, since the distribution $(\|x\|^{-p})^{\wedge}$
is homogeneous of degree $-n+p,$ there exists a distribution
$\gamma$ on the sphere $\Omega$ so that, for every test function
$\phi,$
$$ \langle (\|x\|^{-p})^{\wedge}, \phi \rangle = 
\langle \gamma, \int_{\Bbb R} |t|^{p-1} \phi(t\xi)\ dt \rangle.$$
Applying the latter equality to the test functions of the
form (3), and using the fact  that the
distribution $(\|x\|^{-p})^{\wedge}$ is positive, we conclude that
$\gamma$ is a positive distribution on the sphere.
But every positive distribution is a finite measure on $\Omega,$ 
which follows from an easy argument similar to that in \cite{5, p.143}.
\qed \enddemo

\bigbreak

We need the following simple fact.

\proclaim{Lemma 1} Let $p_k,\ k\in \Bbb N$ be a sequence of
numbers from the interval $(0,n)$ so that the limit
$p= \lim_{k\to \infty} p_k$ exists and $0<p<n.$ 
Suppose that an $n$-dimensional homogeneous
space $B$ embeds in $L_{-p_k}$ for every $k\in \Bbb N.$
Then $B$ embeds in $L_{-p}.$
\endproclaim

\demo{Proof} We can assume that there exists $\epsilon >0$
so that $0< p_k < p+\epsilon <n$ for every $k.$ Fix a non-negative
test function $\phi.$ By Theorem 1,
$\langle (\|x\|^{-p_k})^{\wedge}, \phi \rangle \ge 0.$
Define a function $g$ on $\Bbb R^n$ by 
$g(x)= \|x\|^{-p-\epsilon} |\hat\phi(x)|$ if $\|x\| \le 1,$
and $g(x) = |\hat\phi(x)|$ if $\|x\| > 1.$ The function $g$
is integrable on $\Bbb R^n$ and, for every $k\in \Bbb N,
\ x\in \Bbb R^n,$ we have $g(x)\ge \|x\|^{-p_k}|\hat\phi(x)|.$
By the dominated convergence theorem,
$$\langle (\|x\|^{-p})^{\wedge}, \phi \rangle =
\int_{\Bbb R^n} \|x\|^{-p} \hat\phi(x)\ dx =$$
$$\lim_{k\to \infty} \int_{\Bbb R^n} \|x\|^{-p_k} \hat\phi(x)\ dx =
\lim_{k\to \infty}\langle (\|x\|^{-p_k})^{\wedge}, \phi \rangle \ge 0,$$
and the result follows from Theorem 1.
\qed \enddemo

\bigbreak

In order to prove that every normed space embeds in every 
$L_{-p}$ with $p\in [n-1,n),$ we use the 
following simple facts taken from \cite{11, Lemmas 3,4}.

\proclaim{Lemma 2} Let $p\in (n-1,n)$ and let $f$ 
be an even homogeneous function of degree
$-p$ on $\Bbb R^n\setminus \{0\}$ such that 
$f\vert_\Omega \in L_1(\Omega).$ 
Then for every $\xi\in \Bbb R^n$
$$\hat{f}(\xi) = {{\pi}\over {c}} 
\int_\Omega |(\theta,\xi)|^{-n+p} f(\theta)\ d\theta.$$
where $c = 2^{-n+p+1}\sqrt{\pi} 
\Gamma((-n+p+1)/2)/\Gamma((n-p)/2) > 0.$ 
In particular, $\hat{f}\vert_\Omega \in L_1(\Omega).$ 
\endproclaim

\proclaim{Lemma 3} 
Let $f$ be an even homogeneous function of degree $-n+1$
on $\Bbb R^n\setminus \{0\}$ so that 
$f\vert_\Omega\in L_1(\Omega).$ Then, for every $\xi\in \Omega,$
$$\hat{f}(\xi) = 
\pi \int_{\Omega\cap \{(\theta,\xi)=0\}} f(\theta)\ d\theta.$$
\endproclaim

\bigbreak

If $B=(\Bbb R^n,\|\cdot\|)$ is a homogeneous space
and $p\in [n-1,n),$ the function $f(x)=\|x\|^{-p}$ 
satisfies the conditions of Lemma 2 or Lemma 3.
Therefore, the Fourier transform $(\|x\|^{-p})^\wedge$
is a homogeneous of degree $-n+p,$ positive, locally integrable
in $\Bbb R^n$ function, and, hence, it is a positive  
distribution. By Theorem 1, 

\proclaim{Corollary 1} Every $n$-dimensional homogeneous space 
embeds in $L_{-p}$ for every $p\in [n-1,n).$
\endproclaim

\bigbreak

Note that, in the case $p=n-1,$ the result of Corollary 1 follows
from the case $p\in (n-1,n)$ and Lemma 1, so using
Lemma 3 is not necessary.

\head 3. Embeddings of the spaces $\ell_q^n,\ 0<q\le \infty$ in $L_{-p}.$
\endhead

Let us first show that each of the spaces $L_{-p}$ is large enough
to contain all finite dimensional subspaces of $L_q,\ 0<q\le 2.$

\proclaim{Theorem 2} Every $n$-dimensional subspace of 
$L_q$ with $0<q\le 2$ embeds in $L_{-p}$ for each $p\in (0,n).$
\endproclaim

\demo{Proof} By a well-known result of P.Levy \cite{12},
for every $n$-dimensional subspace $B=(\Bbb R^n,\|\cdot\|)$
of $L_q$ with $0<q\le 2,$ the function $\exp(-\|x\|^q)$ is 
the Fourier transform of a $q$-stable symmetric measure $\mu$ 
on $\Bbb R^n.$ We have
$$\|x\|^{-p} = {q\over {\Gamma(p/q)}}
\int_0^\infty t^{p-1} \exp(-t^q \|x\|^q)\ dt.$$
For every non-negative test function $\phi,$
$$\langle (\|x\|^{-p})^\wedge, \phi \rangle =
\int_{\Bbb R^n} \|x\|^{-p} \hat\phi(x)\ dx =$$
$${{q}\over {\Gamma(p/q)}}\int_0^\infty t^{p-1} dt\int_{\Bbb R^n} 
\hat\phi(x) \exp(-t^q \|x\|^q)\ dx=$$
$${{q}\over {\Gamma(p/q)}}\int_0^\infty t^{p-1} dt\int_{\Bbb R^n} 
\phi(tx) \ d\mu(x) \ge 0.$$
Therefore, $(\|x\|^{-p})^\wedge$ is a positive distribution.
\qed \enddemo
 
\bigbreak
Our next goal is to show that the spaces $\ell_q^n,\ 2<q\le \infty$ 
embed in $L_{-p}$ if and only if $p\in [n-3,n).$ We start with 
calculating the Fourier transform of $\|x\|_\infty^{-p}.$

\proclaim{Lemma 4} If $p\in (0,n)$ then,
for every $\xi\in \Bbb R^n$ with non-zero coordinates, 
$$(\|x\|_\infty^{-p})^{\wedge} (\xi) =
2^{n}p \int_0^\infty t^{-p-1} \prod_{k=1}^n
{{\sin(t\xi_k)}\over {\xi_k}} \ dt. \tag{4}$$

\endproclaim

\demo{Proof} For every $x\in \Bbb R^n,\ x\neq 0,$ we have
$$\|x\|_\infty^{-p} =  
p \int_0^\infty z^{p-1} \chi(z\|x\|_\infty)\ dz,$$
where $\chi$ is the indicator of $[-1,1],$ 
and the integral converges because $p>0.$ Clearly,
$\chi(z\|x\|_\infty) = \prod_{k=1}^n 
\chi(zx_k).$ Therefore, $(\chi(z\|x\|_\infty))^\wedge(\xi) = \prod_{k=1}^n
{{2\sin(\xi_k/z)}\over{\xi_k}}.$ Since $0<p<n,$ for every
test function $\phi\in \Cal S(\Bbb R^n)$ the integral
$$\langle (\|x\|^{-p})^{\wedge}, \phi \rangle = 
p\int_{\Bbb R^n} \hat\phi(x)\ dx 
\int_0^\infty z^{p-1} \chi(z\|x\|_\infty)\ dz\tag{5}$$
converges absolutely, and we can use the Fubini theorem, the
definition of the Fourier transform of distributions and the 
change variables $t=1/z$ to 
show that the expression in the right-hand side of (5) is equal to
$$ 2^{n}p \int_{\Bbb R^n} \phi(\xi)\ d\xi 
\int_0^\infty t^{-p-1} 
\prod_{k=1}^n {{\sin(t\xi_k)}\over {\xi_k}} \ dt,$$
which proves (4). Note that the integral in (4) converges 
absolutely because $0<p<n$ and $\xi$ has non-zero coordinates.
\qed \enddemo

\bigbreak

\proclaim{Lemma 5} If $0<p<n$ then the function $(\|x\|_\infty^{-p})^\wedge$ 
is locally integrable on $\Bbb R^n.$ \endproclaim

\demo{Proof} Since the function $(\|x\|_\infty^{-p})^\wedge$ is homogeneous
of degree $-n+p\in (-n,0),$ it is enough to show that this function
is absolutely integrable on the unit cube $Q_n$ in $\Bbb R^n.$ 
By Lemma 4, 
$$\int_{Q_n} \big| (\|x\|_\infty^{-p})^{\wedge}(\xi)\big|\ d\xi  \le
2^n p \int_0^\infty |t|^{-p-1} 
\Big( \prod_{k=1}^n \int_{-1}^1 \big| {{\sin t\xi_k}\over {\xi_k}}\big|
\ d\xi_k \Big)\ dt =$$
$$2^n p \int_0^\infty |t|^{-p-1} \Big( \int_{-t}^t 
\big|{{\sin u}\over {u}}\big|\ du \Big)^n\ dt < \infty,$$
because $-n-1< -p-1<- 1,$ and $\int_{-t}^t 
\big|{{\sin u}\over {u}}\big|\ du $ is bounded by $2t$ at zero,
and by $2\ln t$ at infinity.
\qed \enddemo

\bigbreak

The integral (4) can easily be calculated if $-1<p<0$
using the representation of 
$\prod_{k=1}^n \sin(t\xi_k)$ as a sign-changing sum
of sins or cosins.
The resulting formula can be extended analytically by $p$ to
all values of $p$ which are not integers. 
We get the following expression:

\proclaim{Lemma 6} Let $p>0,$ $p$ is not an integer,
and $\xi\in \Bbb R^n$ is a vector with non-zero
coordinates. Then, if $n$ is odd we have
$$(\|x\|_\infty^{-p})^{\wedge} (\xi) = $$
$${{(-1)^{{{n-1}\over 2}} 2^{-p} \sqrt{\pi}\ \Gamma({{-p+1}\over 2})}\over
{\xi_1\dots\xi_n\ \Gamma(p/2)}}
\sum_\delta \delta_1\dots\delta_n 
|\delta_1\xi_1+\dots \delta_n\xi_n|^p 
\sgn(\delta_1\xi_1+\dots \delta_n\xi_n).\tag{6}$$
If $n$ is even
$$(\|x\|_\infty^{-p})^{\wedge} (\xi) =
{{(-1)^{{n\over2}+1} 2^{-p} \sqrt{\pi}\ \Gamma((-p+2)/2)}\over
{\xi_1\dots\xi_n\ \Gamma((p+1)/2)}}
\sum_\delta \delta_1\dots\delta_n 
|\delta_1\xi_1+\dots \delta_n\xi_n|^p,\tag{7}$$
where the sum is taken over all changes of signs
$\delta = (\delta_1,\dots,\delta_n),\ \delta_k =\pm 1,
\ k=1,...,n.$
\endproclaim

One can also deduce (6) and (7) from a more general formula in
\cite{9, Lemma 3.3} which allows to calculate the Fourier 
transform of the functions $f(\|x\|_\infty)$ for a large 
class of functions $f$ (note that a multiplier $(-1)^{n-1}$ 
is missing in the formula in \cite{9}; apply that formula
to the functions $f(t)=|t|^p$ with $p>0$ and use analytic 
continuation by $p;$ see Section 2 for the Fourier transform 
of the function $f(t)=|t|^p$).

\medbreak

Let us find the signs of the sums appearing in Lemma 6.

\proclaim{Lemma 7} Let $n>3$ and $0<p<n-3,$ or $n=3$ and $p<0,$ 
where $p$ is not an integer. Then the functions 
$$g_{n,p}(\xi_1,...,\xi_n) = \sum_\delta \delta_1\dots\delta_n 
|\delta_1\xi_1+\dots \delta_n\xi_n|^p 
\sgn(\delta_1\xi_1+\dots \delta_n\xi_n)$$
and 
$$h_{n,p}(\xi_1,...,\xi_n) = \sum_\delta \delta_1\dots\delta_n 
|\delta_1\xi_1+\dots \delta_n\xi_n|^p$$
are sign-changing on
$\Bbb R_{+}^n = \{\xi\in \Bbb R^n:\ \xi_k>0,\ k=1,...,n \}.$
\endproclaim
 
\demo{Proof} First, let $n=3.$ Then, for every $p<0,$
the numbers $g_{3,p}(3,1,1)$ and $h_{3,p}(3,1,1)$ are 
positive, but the numbers $g_{3,p}(1,3,3)$ and 
$h_{3,p}(1,3,3)$ are negative. 

Now, if $n>3$ and $p<n-3,$ put $k=[p]+1.$ Then $n-k \ge 3.$
We have 
$$h_{n,p}(\xi_1,...,\xi_n)= p(p-1)....(p-[p])
\int_{-\xi_1}^{\xi_1}\dots \int_{-\xi_k}^{\xi_k}
\Big( \sum_{\delta_{j}=\pm 1,\ j>k}
\delta_{k+1}...\delta_n $$
$$\big|\sum_{j=1}^k x_j+ 
\sum_{j=k+1}^n\delta_{j}\xi_{j} \big|^{p-[p]-1} 
sgn^k \big(\sum_{j=1}^k x_j+ 
\sum_{j=k+1}^n\delta_{j}\xi_{j}\big) \Big)\ dx_1...dx_k.$$
Since $p-[p]-1 \in (-1,0),$ we can finish the proof using
the argument for $n=3.$ In fact, 
if we put $\xi_{n-2}=3,\ \xi_{n-1}=\xi_n=1$ and make 
$\xi_1,...,\xi_{n-3}$ small enough, the function under the integral
becomes always positive, and the value of the function $h_{n,p}$ 
is positive. However, if we do the same thing with $\xi_{n-2}=1,\ \xi_{n-1}=\xi_n=3,$ the value of $h_{n,p}$ is negative.
A similar argument works for the functions $g_{n,p}\ .$
\qed \enddemo

\proclaim{Lemma 8} Let $n\ge 2$ and $p\in (n-3,n),$ $p$
is not an integer. Put $u_{n,p}=g_{n,p}$ if $n$ is odd,
and $u_{n,p}=h_{n,p}$ if $n$ is even. Then $u_{n,p}$ is 
a positive function on $\Bbb R_{+}^n$  if $p\in (n-2,n),$
and $u_{n,p}$ is a negative function on $\Bbb R_{+}^n$  
if $p\in (n-3,n-2).$ 
\endproclaim

\demo{Proof} We argue by induction. The case $n=2$ is trivial.
Let $\xi\in \Bbb R_{+}^n.$ Without loss of generality we can 
assume that $\xi_1\ge \xi_n.$ Then
$$u_{n,p} (\xi_1,...,\xi_n) = 
p\int_{-\xi_n}^{\xi_n} 
u_{n-1,p-1}(\xi_1 +x, \xi_2,...,\xi_{n-1})\ dx.\tag{8}$$
Since $\xi_1+x \ge 0,$ the result follows from the induction 
hypothesis. \qed \enddemo

\bigbreak

\proclaim{Theorem 3} Let $0<p<n,\ n\ge 3.$ 
The function $\|x\|_\infty^{-p}$ is a 
positive definite distribution if $p\in [n-3,n),$ and it is
not positive definite if $p\in (0,n-3).$ Therefore, the space
$\ell_\infty^n$ embeds in $L_{-p}$ if and only if $p\in [n-3,n).$  
\endproclaim

\demo{Proof} Let $n-3\le p<n.$ By Lemma 5, both functions 
$\|x\|_\infty^{-p}$ and $(\|x\|_\infty^{-p})^\wedge$ are
locally integrable in $\Bbb R^n.$ Also if $n-3<p<n$ and $p$
is not an integer then, by Lemmas 6 and 8,
the function $(\|x\|_\infty^{-p})^\wedge$ (which is even by
each variable) is non-negative almost everywhere (with respect
to Lebesque measure) on $\Bbb R^n.$ It is easily seen now that, 
for every $p\in (n-3,n)$ which is not an integer,
$(\|x\|_\infty^{-p})^\wedge$ is a positive distribution. 
By Lemma 1, the same is true for $p=n-3, n-2, n-1.$

Let $0<p<n-3.$ By Lemma 7, if $p$ is not an integer then
the function $(\|x\|_\infty^{-p})^\wedge$ has opposite signs
at two different points, and the function is continuous
in neighborhoods of those points, so $(\|x\|_\infty^{-p})^\wedge$ 
is not a positive distribution. We can show the same thing using
a different argument which also applies to the integers $p.$
In fact, if for some $0<p<n-3$ the function 
$(\|x\|_\infty^{-p})^\wedge$ is non-negative almost everywhere, 
then by Lemma 11 so is the function $(\|x\|_q^{-p})^\wedge$
which contradicts Lemma 10.
\qed \enddemo

\bigbreak

Let us pass to the spaces $\ell_q^n,\ 2<q<\infty.$
Denote by $\gamma_q$ the Fourier transform of
the function $z\to \exp(-|z|^q),\ z\in \Bbb R.$
The properties of the functions $\gamma_q$ were studied
by Polya \cite{16}. In particular, if $q$ is not an even integer,
the function $\gamma_q(t)$
behaves at infinity like $|t|^{-q-1}.$ Namely 
(see \cite{17, Part 3, Problem 154}),
$$\lim\limits_{t\to\infty} t^{1+q}\gamma_q(t)=
2\Gamma (q+1) \sin (\pi q/2).$$
If $q$ is an even integer, the function $\gamma_q$
decreases exponentially at infinity.
The integral 
$$S_q(\alpha)=\int_{\Bbb R}|t|^\alpha \gamma_q(t)\ dt $$
converges absolutely for every $\alpha\in (-1,q).$ These 
moments can easily be calculated (see \cite{20} or
\cite{7}; $\alpha$ is not an even integer):
$$S_q(\alpha)=2^{\alpha+2}\pi^{1/2} \Gamma(-\alpha / q)
\Gamma ((\alpha+1)/2)/ (q\Gamma(-\alpha/2)).$$
Clearly, the moment $S_q(\alpha)$ is positive if 
$\alpha\in (-1,0)\cup (0,2),$ and the moment is negative 
if $\alpha\in (2,\min(q,4)).$

The Fourier transform of the function $\|x\|_q^\beta$
was calculated in \cite{7}.
\proclaim{Lemma 9}
Let $q>0$, $n\in\Bbb N$, $-n<\beta<qn$, 
$\beta/q \not\in \Bbb N\cup \{0\}$,
$\xi=(\xi_1, \ldots, \xi_n)\in \Bbb R^n$, 
$\xi_k\not= 0$, $1\leq k\leq n$.
Then
$$ (\|x\|_q^\beta)^{\wedge}(\xi) =
\frac{q}{\Gamma(-\beta/q)} 
\int_0^\infty t^{n+\beta-1}\prod_{k=1}^n
\gamma_q (t\xi_k)\, dt.
$$
\endproclaim
\bigbreak

Let us prove that the function $(\|x\|_q^{-p})^\wedge$ changes 
its sign if $0<p<n-3.$ The following argument is similar to 
that used in the proof of the 1938 Schoenberg's conjecture 
on positive definite functions in \cite{7}.

\proclaim{Lemma 10} If $q>2,\ n>3,\ p\in (0,n-3)$ 
then the distribution $\|x\|_q^{-p}$ is not positive definite.
\endproclaim

\demo{Proof}  By Lemma 9 and properties of the moments 
$S_q(\alpha),$ the integral
$$I(\alpha_1,\dots,\alpha_{n-1})= \int_{\Bbb R} 
|\xi_1|^{\alpha_1}\dots|\xi_{n-1}|^{\alpha_{n-1}} 
(\|x\|_q^{-p})^\wedge(\xi_1,...,\xi_{n-1},1)\ d\xi_1\dots d\xi_{n-1}=$$
$$ S_q(\alpha_1)\dots S_q(\alpha_{n-1}) 
S_q(-\alpha_1-\dots-\alpha_{n-1}-p)$$
converges absolutely if the numbers 
$\alpha_1,\dots,\alpha_{n-1},-\alpha_1-\dots-\alpha_{n-1}-p$
belong to the interval $(-1,q).$ Choosing $\alpha_k\in (-1,0)$ for 
every $k=1,...,n-1,$ we have the moments $S_q(\alpha_k),\ k=1,...,n$ 
positive, and we can make $-\alpha_1-\dots-\alpha_{n-1}-p$ equal to any
number from $(-p,n-1-p)\cap (-1,q).$ This interval contains a
neighborhood of 2, and, since the moment function 
$S_q$ changes its sign at 2, we can make the integral
$I(\alpha_1,\dots,\alpha_{n-1})$ positive for one
choice of $\alpha$'s and negative for another choice.
This means that the function $(\|x\|_q^{-p})^\wedge$
is sign-changing.
\qed \enddemo

\bigbreak
To show that, for $p\in [n-3,n),$ the function 
$(\|x\|_q^{-p})^\wedge$ is positive almost everywhere,
we first express this function in terms of the function 
$(\|x\|_\infty^{-p})^\wedge,$ and then positivity will
follow from Lemma 8 (for those $p$ which are not integers).

\proclaim{Lemma 11} Let $q>0,\ p\in (0,n).$ Then, for every
$\xi\in \Bbb R^n$ with non-zero coordinates,
$$(\|x\|_q^{-p})^\wedge (\xi) =
{{q^{n+1}}\over {p\Gamma(p/q)}}\int_0^\infty \dots \int_0^\infty$$ 
$$(t_1\dots t_n)^q \exp(-\|t\|_q^q) 
(\|x\|_\infty^{-p})^\wedge(t_1\xi_1,...,t_n\xi_n)\ dt_1...dt_n.$$
\endproclaim

\demo{Proof} For every $x\in \Bbb R,$ we have
$$\exp(-|x|^q) = q\int_0^\infty \chi(ux)
\ u^{-1-q}\exp(-u^{-q})\ du,\tag{9}$$
where, as before, $\chi$ is the indicator of $[-1,1].$

The Fourier transform of the function $x\mapsto \chi(ux)$
is equal to $(\chi(ux))^\wedge(\xi) = 2\sin(\xi/u)/\xi.$
Calculating the Fourier transforms of both sides of (9)
and making the change of variables $t=1/u,$ 
we get an integral representation for the function $\gamma_q:$ 
for every $\xi\in \Bbb R,$
$$\gamma_q(\xi)= 2q\int_0^\infty {{\sin(t\xi)}\over {t\xi}} 
t^q \exp(-t^q)\ dt.$$
By Lemma 9,
$$ (\|x\|_q^{-p})^{\wedge}(\xi) =
\frac{q}{\Gamma(p/q)} 
\int_0^\infty z^{n-p-1}\prod_{k=1}^n
\gamma_q (z\xi_k)\, dz = $$
$${{2^n q^{n+1}}\over {\Gamma(p/q)}} \int_0^\infty\dots \int_0^\infty  
(t_1\dots t_n)^q \exp(-\|t\|_q^q)
\int_0^\infty z^{n-p-1} \prod_{k=1}^n 
{{\sin(t_k\xi_k z)}\over {t_k\xi_k z}}\ dz,$$
and the result follows from Lemma 4.
\qed \enddemo
\bigbreak

\proclaim{Theorem 4} Let $0<p<n,\ n\ge 3.$ If $2<q< \infty$ then  
$\|x\|_q^{-p}$ is a positive definite distribution 
if $p\in [n-3,n),$ and it is
not positive definite if $p\in (0,n-3).$ Therefore, the space
$\ell_q^n$ embeds in $L_{-p}$ if and only if $p\in [n-3,n).$  
\endproclaim

\demo{Proof} If $p\in (n-3,n)$ and $p$ is not an integer,
then, by Lemmas 11 and 8 the function $(\|x\|_q^{-p})^\wedge$
is positive almost everywhere. It is also locally integrable
which follows from Lemma 9 and an argument similar to that in Lemma 5.
Therefore, $(\|x\|_q^{-p})^{\wedge}$ is a positive distribution,
and, by Theorem 1, the space $\ell_q^n$ embeds in $L_{-p}.$ 
One can use Lemma 1 to add the integers $n-3, n-2$ and $n-1.$ 

In the case where $0<p<n-3,$ the result follows from Lemma 10.
\qed \enddemo

\head 4. Inequalities of correlation type for the expectations
of norms of stable vectors. \endhead

For $0<q\le 2,$ let $X=(X_1,...,X_n)$ 
be a symmetric $q$-stable random vector
which means that 
the characteristic functional of the vector $X$ has the form
$$\phi(\xi)= \exp(-\|\sum_{i=1}^n \xi_is_i\|_q^q),\quad 
\xi\in \Bbb R^n, \tag{10}$$
where $s_1,\dots,s_n\in L_q([0,1]).$ 
In this section, we use the notation $\|\cdot\|_q$ 
for the norm of the space $L_q([0,1]).$

Fix an integer $k,\ 1\le k < n,$ and consider the set 
$\Cal A(X,k)$ of all $n$-dimensional symmetric $q$-stable 
random vectors whose first $k$ coordinates have the same 
joint distribution as $X_1,...,X_k,$ and whose last $n-k$
coordinates have the same joint distribution as 
$X_{k+1},...,X_n.$ We denote by $Y=(Y_1,...,Y_n)$ the 
vector from $\Cal A(X,k)$ for which every $Y_i$ and $Y_j$
with $1\le i \le k,\ k+1\le j \le n$ are independent.
Then, the characteristic functional of $Y$ is equal to
$$\phi_0(\xi)= \exp(-\|\sum_{i=1}^k \xi_is_i\|_q^q - 
\|\sum_{i=k+1}^n \xi_is_i\|_q^q).$$

\bigbreak

Given an $n$-dimensional homogeneous space
$B=(\Bbb R^n,\|\cdot\|)$ and a real number $p,$ 
we are interested in conditions on $B$ and $p$ under 
which the independent case is extremal in the sense that 
the expectation $\Bbb E(\|Y\|^p)$ is the minimal or maximal 
value of $\Bbb E(\|Z\|^p),\ Z\in \Cal A(X,k).$

\bigbreak

First, we consider the 
case where $B=(\Bbb R^n,\|\cdot\|)$ is an $n$-dimensional
subspace of $L_p,\ p>0$ satisfying the following symmetry condition: 
for every $u\in \Bbb R^k,\ v\in \Bbb R^{n-k},$
$$\|(u,v)\| = \|(u,-v)\|.
\tag{*}$$
We use the representation (1) and a standard argument from 
the theory of stable processes to show that, if $0<p\le q$ then 
$$\Bbb E(\|Y\|^p)= \max\{\Bbb E(\|Z\|^p):\ Z\in \Cal A(X,k)\}.$$

As it was mentioned in the Introduction, the condition that 
$B$ is a subspace of $L_p$ is restricting, for example, 
the most interesting case of $B=\ell_\infty^n$ is not covered.
However, we replace the standard argument by the technique of 
embedding in $L_{-p}$, which allows to get
more spaces involved. We prove that if $B$ embeds into 
$L_{-p},\ p\in (0,n)$ and has the symmetry (*) then
$$\Bbb E(\|Y\|^{-p})= \min\{\Bbb E(\|Z\|^{-p}):\ Z\in \Cal A(X,k)\}.$$

\bigbreak

Let us start with the standard technique.
If $B$ is a subspace of $L_p$ with $p>0,$
then one can use the well-known formula for the expectations
of the scalar products of $q$-stable vectors with fixed vectors
to reduce the estimation of $\Bbb E(\|X\|^p)$ to simple 
properties of the $L_q$-norms.

We need a few simple inequalities for the $L_q$-norms which follow
from Clarkson's inequality (see \cite{2}). For the reader's 
convenience we include the proof. 

\proclaim{Lemma 12} Let $x,y\in L_q([0,1]),\ 0<q\le 2.$ Then
$$\exp(-\|x+y\|_q^q) + \exp(-\|x-y\|_q^q) \ge
2\exp(-\|x\|_q^q-\|y\|_q^q). \tag{11}$$
Also for every $0<p\le q$
$$\|x + y\|_q^p +
\|x - y\|_q^p \le
2(\|x\|_q^q + \|y\|_q^q)^{p/q}.\tag{12}$$
Finally, for $q=2$ and $p>2$ the inequality (12) reverses.
\endproclaim

\demo{Proof} First, note that for any $0<q\le 2$
$$\|x+y\|_q^q + \|x-y\|_q^q \le 2(\|x\|_q^q + \|y\|_q^q), \tag{13}$$
and this is a simple consequence of the same inequality for real
numbers. Now to get (11) apply the relation between the
arithmetic and geometric means and then use (13). The inequality (12)
also follows from (13):
$$ \Big( {{\|x + y\|_q^p + \|x - y\|_q^p}\over 2} \Big)^{1/p} \le
\Big( {{\|x + y\|_q^q + \|x - y\|_q^q}\over 2} \Big)^{1/q} \le
(\|x\|_q^q + \|y\|_q^q)^{1/q}.$$
Finally, if $q=2$ the latter calculation works for $p>2$ where
the first inequality goes in the opposite direction, and the 
second inequality turns into an equality.
\qed \enddemo

\bigbreak

\proclaim{Proposition 1} Let $q,k,X,Y$ be as in the beginning 
of this section, $0<p\le q.$ Let $B=(\Bbb R^n,\|\cdot\|)$ be a 
subspace of $L_p$ satisfying the condition (*).
Then 
$$\Bbb E\ (\|Y\|^p)=\max\{\Bbb E\ (\|Z\|^p):\ Z\in \Cal A(X,k)\}.$$ 
Also, if $q=2$ and $p>2$ then $\Bbb E\ (\|Y\|^p)$ is the minimal value. 
\endproclaim
 
\demo{Proof} A basic property of the stable vector with the characteristic 
function (10) is that, for any
vector $\xi\in \Bbb R^n,$ the random variable $(X,\xi)$ has the same
distribution as $\|\sum_{i=1}^n \xi_is_i\|_q U,$ where $U$ is the
standard one-dimensional $q$-stable random variable.
Therefore, if $p<q$  then 
$$\Bbb E\ |(X,\xi)|^p = c_{p,q}\|\sum_{i=1}^n \xi_is_i\|_q^p,\tag{14}$$
where $c_{p,q}$ is the $p$-th moment of $U$
(which exists only for $p<q$ if $q<2,$ and
it exists for every $p>0$ if $q=2;$ see \cite {20} for
a formula for $c_{p,q}).$ Similarly, we get
$$\Bbb E\ |(X_{-},\xi)|^p = c_{p,q}\|\sum_{i=1}^k \xi_is_i- 
\sum_{i=k+1}^n \xi_is_i\|_q^p, $$ 
where $X_{-}=(X_1,...,X_k,-X_{k+1},...,-X_n).$ Also,
$$\Bbb E\ |(Y,\xi)|^p = c_{p,q}(\|\sum_{i=1}^k \xi_is_i\|_q^q + 
\|\sum_{i=k+1}^n \xi_is_i\|_q^q)^{p/q}.$$

Since $(\Bbb R^n,\|\cdot\|)$ is a subspace of $L_p([0,1]),$
we can use the Blaschke-Levy representation (1) and after that 
the formula (14) to get

$$\Bbb E\ (\|X\|^p) = 
\int_S \Bbb E\ (|(X,\xi)|^p)\ d\gamma(\xi)=
c_{p,q}\int_S \|\sum_{i=1}^n \xi_is_i\|_q^p\ d\gamma(\xi).\tag{15}$$
Similarly,
$$E (\|Y\|^p)= c_{p,q}\int_S (\|\sum_{i=1}^k \xi_is_i\|_q^q + 
\|\sum_{i=k+1}^n \xi_is_i\|_q^q)^{p/q}\ d\gamma(\xi),\tag{16}$$
$$\Bbb E\ (\|X_{-}\|^p)= c_{p,q}\int_S \|\sum_{i=1}^k \xi_is_i- 
\sum_{i=k+1}^n \xi_is_i\|_q^p \ d\gamma(\xi).\tag{17}$$
Since $0<p\le q,$ the equalities (15), (16), (17) in conjunction 
with (12) imply 
$\Bbb E\ (\|X\|^p)+\Bbb E\ (\|X_{-}\|^p)\le 2\Bbb E\ (\|Y\|^p),$
and now the result follows from the property of
the norm that $\|X\|=\|X_{-}\|.$ In the case $q=2,\ p>2$ we use
the corresponding part of Lemma 12.
\qed \enddemo

\bigbreak

\subheading{Remark}  For $p>q,\ q< 2 $ the expectation
of $\|X\|^p$ does not exist so the statement of Proposition 1
does not make sense in that case.

\bigbreak

\proclaim{Theorem 5} Let $q,k, X,Y$ be as in Proposition 1,
and suppose that $0<p<n$ and $B=(\Bbb R^n,\|\cdot\|)$ is a homogeneous
space which embeds in $L_{-p}$ and whose norm satisfies the 
symmetry condition (*). Then
$\Bbb E (\|X\|^{-p}) \ge \Bbb E (\|Y\|^{-p}).$ \endproclaim

\demo{Proof} By Theorem 1, the function $\|x\|^{-p}$
is a positive definite distribution, and by a generalization
of Bochner's theorem \cite{5}, this function is the Fourier transform 
of a tempered measure $\mu$ on $\Bbb R^n.$

Let $P_X$ be the $q$-stable measure in $\Bbb R^n$ 
according to which the random vector $X$ is distributed.
Applying the Parseval equality and formula (10) for the 
characteristic function of $X$ we get
$$\Bbb E (\|X\|^{-p})= \int_{\Bbb R^n} \|x\|^{-p}\ dP_X(x) =
\int_{\Bbb R^n} \widehat{P_X}(\xi)\ d\mu(\xi) =$$
$$\int_{\Bbb R^n} \ exp(-\|\sum_{i=1}^n \xi_is_i\|_q^q)\ d\mu(\xi).$$
Note that the function $\|x\|^{-p}$ is locally integrable in $\Bbb R^n$ 
because $0<p<n.$ Similarly,
$$\Bbb E (\|X_{-}\|^{-p})= 
\int_{\Bbb R^n} \ exp(-\|\sum_{i=1}^k \xi_is_i-
\sum_{i=k+1}^n \xi_is_i\|_q^q)\ d\mu(\xi),$$
where $X_{-}=(X_1,...,X_k,-X_{k+1},...,-X_n),$ and 
$$\Bbb E (\|Y\|^{-p})= 
\int_{\Bbb R^n} \ exp(-\|\sum_{i=1}^k \xi_is_i\|_q^q-
\|\sum_{i=k+1}^n \xi_is_i\|_q^q)\ d\mu(\xi).$$
Now by the inequality (11) from Lemma 12 and taking in account that
$\mu$ is a positive measure, we get
$$\Bbb E (\|X\|^{-p}) + \Bbb E (\|X_{-}\|^{-p}) \ge 
2\Bbb E (\|Y\|^{-p}),$$
and the result follows from the property (*).\qed \enddemo

\bigbreak

The following is an immediate consequence of Theorem 5 in 
conjunction with Theorems 2,3,4 and Corollary 1.

\proclaim{Corollary 2} Let $B=(\Bbb R^n,\|\cdot\|)$ be a homogeneous
space, $0<p<n,$ and $q,k,X,Y$ as above. Then the inequality
$$\Bbb E (\|X\|^{-p}) \ge \Bbb E (\|Y\|^{-p})$$
holds in each of the following cases:
\item{(i)} $B$ is any $n$-dimensional homogeneous space
satisfying the condition (*) and $p\in [n-1,n);$
\item{(ii)} $B$ is an $n$-dimensional subspace of $L_r$
with $0<r\le 2$ satisfying the condition (*) and $p$ is 
any number from $(0,n);$
\item{(iii)} $B=\ell_q^n,\ n\ge 3,\ 2<q\le \infty$
and $p\in [n-3,n).$
\endproclaim

\subheading{Acknowledgements} Part of this work was done
during the NSF Workshop on Linear Analysis and Probability
held at Texas A\&M University in 1995 and 1996. I would like 
to thank G. Schechtman, T. Schlumprecht and J. Zinn for very 
useful discussions.


\Refs 

\ref \no 1\paper Lois stables et espaces $L_p$ \by J. Bretagnolle,
D. Dacunha-Castelle and J. L. Krivine
\jour Ann. Inst. H. Poincar\'e  Probab.  Statist. \vol 2 \yr 1966 
\pages 231--259 \endref

\ref \no 2 \by J.A. Clarkson \paper Uniformly convex spaces \yr 1936
\jour Trans. Amer. Math. Soc. \vol 40  \pages 396--414 \endref

\ref \no 3\paper A representation of the symmetric bivariate 
Cauchy distributions
\by T. S. Ferguson \jour Ann. Math. Stat. 
\vol 33 \yr 1962 \pages 1256--1266 \endref  

\ref \no 4 \by I. M. Gelfand and G. E. Shilov \book Generalized functions 1.
Properties and operations  \publ Academic Press \publaddr New York \yr 1964 \endref

\ref \no 5 \by I. M. Gelfand and N. Ya. Vilenkin 
\book Generalized functions 4. Applications of harmonic analysis 
\publ Academic Press \publaddr New York \yr 1964 \endref

\ref \no 6 \paper A class of negative definite functions
\by C. Herz \jour Proc. Amer. Math. Soc. 
\vol 14 \yr 1963 \pages 670--676 \endref 

\ref \no 7 \paper Schoenberg's problem on positive definite functions
\by A. Koldobsky
\jour Algebra and Analysis \vol 3 \yr 1991 \pages 78--85 
\paperinfo ( English translation in
St. Petersburg Math. J. 3 (1992), 563-570) \endref

\ref \no 8\paper  Generalized Levy representation of norms 
and isometric embeddings into $L_p$-spaces \by A. Koldobsky
\jour  Ann. Inst. H.Poincare ser.B \vol 28 \yr 1992 \pages 335--353 
\endref

\ref \no 9 \paper Characterization of measures by potentials
\by A. Koldobsky
\jour J. Theor. Prob. \vol 7 \yr 1994 \pages 135-145 \endref

\ref \no 10 \by A. Koldobsky 
\paper Positive definite functions, stable measures, and
isometries on Banach spaces 
\jour Lect. Notes in Pure and Appl. Math. \vol 175 
\yr 1995 \pages 275--290 \endref

\ref \no 11 \by A. Koldobsky
\paper Inverse formula for the Blaschke-Levy representation
and its applications to zonoids and sections of star bodies
\jour  file  koldobskyinvzonsect.tex 
on the Banach Space Bulletin Board
\endref

\ref \no 12\book Th$\acute {e}$orie de l'addition de variable 
al$\acute {e}$atoires
\by P. Levy \publ Gauthier-Villars \publaddr Paris \yr 1937  
\endref  

\ref \no 13\paper On the extension of operators with finite 
dimensional range \by J. Lindenstrauss
\jour Illinois J. Math. \vol 8 \yr 1964 
\pages 488--499 \endref

\ref \no 14\paper Positive definite functions on $\ell_{\infty}$
\by J. Misiewicz \jour Statist. Probab. Lett. \vol 8 \yr 1989 
\pages 255--260 \endref  

\ref \no 15 \by J. Misiewicz
\paper Sub-stable and pseudo-isotropic processes,
\ preprint  \endref

\ref \no 16\paper On the zeroes of an integral function represented by Fourier's integral \by G. Polya \yr 1923\jour Messenger Math. \vol 52  \pages 185--188  \endref  

\ref \no 17\by G. Polya and G. Szego \yr 1964 
\book Aufgaben und lehrsatze aus der analysis 
\publ Springer-Verlag \publaddr Berlin-New York \endref

\ref \no 18 \by G. Schechtman, T. Schlumprecht  and J. Zinn
\paper On the Gaussian measure of the intersection
of symmetric convex sets,\ preprint  \endref

\ref \no 19\paper Metric spaces and positive definite functions
\by I. J. Schoenberg \yr 1938\jour Trans. Amer. Math. Soc. 
\vol 44 \pages 522--536 \endref  

\ref \no 20 \book One-dimensional stable distributions
\by V. M. Zolotarev \publ Amer. Math. Soc. 
\publaddr Providence \yr 1986 
\endref  















\endRefs
\enddocument